\documentclass{article}
\usepackage{amssymb, amsmath}
\usepackage{amsthm}
\usepackage{epsfig}

\def\v{{\varphi}}

\begin{document}

\title{ \bf An example for nonequivalence of
symplectic capacities}
\author{Ursula Hamenst\"adt \thanks
{Research partially supported by SFB 256 and SFB 611.}
\\
Mathematisches Institut der Unversit\"at Bonn\\
Beringstra\ss{}e 1, D-53115 Bonn, Germany\\
e-mail: ursula@math.uni-bonn.de
}
\date{}
\maketitle

\begin{abstract} We 
construct an open bounded star-shaped
set
$\Omega\subset \mathbb{R}^{4}$
whose cylindrical capacity is strictly bigger than its
proper displacement energy.
\end{abstract}

\section{Introduction}

Consider the standard $2n-$dimensional euclidean space  
$\mathbb{R}^{2n}$ equipped 
with the euclidean symplectic form $\omega_0 =
\sum^n_{i =1} d x_{2i-1} \wedge d x_{2i}$. In this paper we are
interested in symplectic invariants of nonempty open subsets of
$(\mathbb{R}^{2n},\omega_0)$. One example of such an invariant is a
{\sl relative} or {\sl nonintrinsic capacity} [MS] which
associates to every open subset  $\Omega$ of $\mathbb{R}^{2n}$ a
number $c(\Omega)\in [0,\infty]$. This number $c(\Omega)$ measures
the symplectic size of $\Omega$ in such a way that the following
three properties hold.

\begin{enumerate}
\item[A1] {\sl Monotonicity:} $c(\Omega)\leq c
(D)$ if there is a {\sl global}
symplectomorphism of $\mathbb{R}^{2n}$ which maps
$\Omega$ into $D$.
\item[A2] {\sl Conformality:} $c(a\Omega)=a^2  c(\Omega)$
for all $a>0$.
\item[A3] {\sl Nontriviality:} $c(B^{2n}(1))=1=
c(Z^{2n}(1))$
for the open normalized ball
$B^{2n} (1)$ of radius $\sqrt{1/\pi}$
and the open symplectic cylinder $Z^{2n} (1)=
B^2(1)\times \mathbb{R}^{2n-2}$
in the standard space $(\mathbb{R}^{2n}, \omega_0)$.
\end{enumerate}

Here we use coordinates $(x_1,\dots,x_{2n})$ in $\mathbb{R}^{2n}$ and
we write $B^{2n} (r) = \{ x \in \mathbb{R}^{2n} \mid \vert x\vert^2 <
r/\pi\}$ and $Z^{2n} (r) = B^2(r)\times \mathbb{R}^{2n-2}= \{ x\in
\mathbb{R}^{2n}\mid x_1^2 + x_2^2 < r/\pi \}$ for the ball and
cylinder of capacity $r>0$ in $\mathbb{R}^{2n}$.

The celebrated {\sl non-squeezing lemma} of Gromov [G] shows that
for $r>1$ the ball $B^{2n}(r)$ does not admit a symplectic
embedding into the cylinder $Z^{2n}(1)$. This implies that
relative capacities do exist, and in fact there are many ways to
define them. The resulting invariants do not coincide in general.
We will consider the following four examples of such relative
capacities.

The {\sl Gromov width} assigns to an open set 
$\Omega\subset \mathbb{R}^{2n}$
the supremum $c_0(\Omega)$ 
of all numbers $r>0$ such that there is
a symplectic embedding of the ball $B^{2n}(r)$ into
$\Omega$. By monotonicity, the Gromov width is the smallest
capacity which means that if $c^\prime$ is any relative
capacity, then $c_0(\Omega)\leq c^\prime(\Omega)$ for
every open set $\Omega\subset \mathbb{R}^{2n}$.

Let ${\cal O}$ be the family of nonempty open bounded subsets of
$\mathbb{R}^{2n}$. Our second example is the {\sl cylindrical
capacity} which associates to $\Omega\in {\cal O}$ the infimum
$c_p(\Omega)$ of all numbers $r>0$ for which there is a
symplectomorphism of $\mathbb{R}^{2n}$ which maps $\Omega$ into the
cylinder $Z^{2n}(r)$ [P]. 
If $\Omega\subset \mathbb{R}^{2n}$ is unbounded
then we define $c_p(\Omega)=\sup\{c_p(\Omega^\prime)\mid
\Omega^\prime \in {\cal O}, \Omega^\prime\subset \Omega\}$. By
monotonicity, the cylindrical capacity is the biggest relative
capacity which means that if $c^\prime$ is any relative capacity,
then $c^\prime(\Omega)\leq c_p(\Omega)$ for every $\Omega\in {\cal
O}$.

Third the {\sl displacement energy} is defined as follows. Recall
that a compactly supported smooth time dependent function $H(t,x)$
on $[0,1]\times \mathbb{R}^{2n}$ induces a time-dependent Hamiltonian
flow on $\mathbb{R}^{2n}$. Its time-one map $\v$ is then a
symplectomorphism of $\mathbb{R}^{2n}$. The group ${\cal D}$ of
compactly supported symplectomorphisms obtained in this way is
called the group of {\sl compactly supported Hamiltonians} [HZ].

The {\sl Hofer-norm} on the group ${\cal D}$ assigns to $\v\in
{\cal D}$ the value 
\[\Vert \v\Vert= \inf_H\bigl(\sup_{t\in
[0,1]}(\sup_{x\in \mathbb{R}^{2n}} H(t,x)-\inf_{x\in \mathbb{R}^{2n}}
H(t,x))\bigr)\] 
where $H$ ranges over the set of all compactly
supported time dependent functions whose Hamiltonian flows induce
$\v$ as their time-one map. The Hofer-norm $\Vert\; \Vert$ induces a
bi-invariant distance function $d$ on the group ${\cal D}$ by
defining $d(\v,\psi)= \Vert \v\circ \psi^{-1}\Vert$, in particular we
have $\Vert \varphi \Vert > 0$ for $\varphi \ne Id$ [HZ].

For a bounded set $A\subset \mathbb{R}^{2n}$ we define the
{\sl displacement energy}
$d(A)$ to be the infimum of the Hofer
norms $\Vert \v\Vert$ of all those $\v\in {\cal D}$ which
{\sl displace} $A$,
i.e. for which we
have $\varphi(A)\cap A=\emptyset$.
If $A\subset \mathbb{R}^{2n}$ is unbounded then we
define $d(A)=\sup\{d(A^\prime)\mid A^\prime\subset A,
A^\prime \,{\rm bounded}\}.$

Since the cube $(0, 1) \times (0, a) \subset\mathbb{R}^{2}$ of area
$a > 0$ is displaced by the time-one map of the Hamiltonian flow
induced by the time-independent function $H (t,x, y) = ax$, the
displacement energy of the cylinder $(0, 1) \times (0,
a)\times\mathbb{R}^{2n-2} \subset\mathbb{R}^{2n}$ in $\mathbb{R}^{2n}$  is
not bigger than its capacity $a > 0$. This implies in particular
that $d(\Omega)\leq c_p(\Omega)$ for every open bounded subset of
$\mathbb{R}^{2n}$. On the other hand, the displacement energy of an
euclidean ball of capacity $a$ is not smaller than $a$ (this was
first shown by Hofer; we refer to [HZ] and [LM] for proofs and
references). Since moreover clearly $d (\Omega^\prime) \le d
(\Omega)$ if $ \Omega^\prime\subset \Omega$, the displacement
energy is a relative capacity.

Following [HZ] we call two
subsets $A, B$ of $\mathbb{R}^{2n}$ {\sl properly separated} if there
is a symplectomorphism
$\Psi$ of $\mathbb{R}^{2n}$ such that $\overline{\Psi (A)}
\subset \{x_1 < 0\}$ and $\overline{\Psi (B)} \subset \{x_1 >
0\}$. Define the {\sl proper displacement energy} $e (\Omega)$
of a set $\Omega \in {\cal O}$ to be the infimum of the
Hofer-norms $\Vert \varphi \Vert$ of all those $\varphi \in {\cal D}$
for which $\varphi (\Omega)$ and $\Omega$ are properly separated.
If $\Omega \subset \mathbb{R}^{2n}$ is unbounded we define $e
(\Omega) = \sup \{ e (\Omega^\prime) \mid \Omega^\prime\subset
\Omega,\Omega^\prime\; {\rm bounded } \}.$ As before, the
proper displacement energy is a relative capacity. We have
the inequalities $c_0(\Omega)\leq d(\Omega)\leq e(\Omega)\leq
c_p(\Omega)$ for every set $\Omega\in {\cal O}$.

Even for star-shaped subsets of $\mathbb{R}^{2n}$
$(n\geq 2)$  
our above capacities define different symplectic invariants.
The earliest result known to me in this direction is
due to Hermann [He].
He constructed
for every $n\geq 2$ star-shaped Reinhardt-domains in
$\mathbb{R}^{2n}$ with arbitrarily small volume and hence
arbitrarily small Gromov width whose displacement energy
is bounded from below by $1$. For the 
estimate of the displacement energy he
uses a remarkable
result of Chekanov [C] who showed that the displacement
energy of a closed Lagrangian submanifold of $\mathbb{R}^{2n}$ is
positive. 

The displacement energy of closed Lagrangian submanifolds is
not the only obstruction for embeddings of a star-shaped set
$\Omega$ into a cylinder of small capacity. The purpose
of this note is to show.

\bigskip

{\bf Theorem :} {\it There is an
open bounded starshaped subset $\Omega$ of $\mathbb{R}^{4}$
with $e(\Omega)<c_p(\Omega).$}

\bigskip

A modification of our construction can be used to obtain
for every $n\geq 2$ 
examples of open bounded subsets $\Omega$ 
of $\mathbb{R}^{2n}$ with 
$e(\Omega)<c_p(\Omega)$.

The organization of this note is as follows.
In Section 2 we collect some results on symplectic
embeddings and symplectic isotopies which are needed
for the proof of our 
Theorem. The theorem is proved in Section 3, and
Section 4 contains some additional remarks 
on our capacities.

\section{Extensions of symplectic embeddings and isotopies}

In this section we formulate some versions of well known
existence results for
extensions of symplectic embeddings which
are needed for the construction of our example.

Denote again by ${\cal O}$ the collection of 
all open bounded subsets of $\mathbb{R}^{2n}$.
Define a {\sl proper symplectic embedding} of 
an open bounded set $\Omega\in {\cal O}$ into
an open (not necessarily bounded) subset 
$C$ of $\mathbb{R}^{2n}$ 
to be a symplectic embedding of a neighborhood of $\Omega$ 
in $\mathbb{R}^{2n}$ into $C$. 

The {\sl neighborhood extension theorem}
of Banyaga [B] (see also [MS] for a proof) gives a
sufficient condition for the existence of a symplectomorphism
of $\mathbb{R}^{2n}$ extending a given proper symplectic embedding
of a suitable set $\Omega\in {\cal O}$.

\bigskip
{\bf Theorem 2.1:} {\it 
Let $\Omega\in {\cal O}$ be 
an open bounded set such that $H^1(\overline{\Omega},\mathbb{R})=0.$
Then for every symplectic embedding $\v$ of a neighborhood of
$\overline{\Omega}$ into $\mathbb{R}^{2n}$ there is a symplectomorphism
of $\mathbb{R}^{2n}$ which coincides with $\v$ near $\overline{\Omega}$.}

\bigskip

Let again $\Omega\in {\cal O}$. 
A {\sl strict symplectomorphism} of $\Omega$ is a symplectomorphism
$\v$ of $\Omega$ which equals the identity near the boundary
of $\Omega$.
A {\sl strict isotopy} of $\Omega$ is
a $1$-parameter family $\v_t$ of symplectomorphisms of $\Omega$
which coincides with the identity near the boundary
of $\Omega$ and such that $\v_0=Id$. 
We did not find the precise formulation of the following 
lemma
in the literature, so we include the easy proof for convenience.

\bigskip

{\bf Lemma 2.2:} {\it Let $\Omega\subset \mathbb{R}^{2n}$ be 
open, bounded and star-shaped. Then any two strict symplectomorphisms
of $\Omega$ are strictly isotopic.}

{\it Proof:} 
Let $\Omega\subset \mathbb{R}^{2n}$ be open, bounded and star-shaped
with respect to the origin. 
We have to show that every strict symplectomorphism
$\Psi$ of $\Omega$ is strictly isotopic to the identity.

For this we use the arguments of Banyaga. Namely, since the
support of $\Psi$ is compact there is an isotopy
$\Psi_t$ of the identity with compact support in $a\Omega$ for some
$a\geq 1$ and such that $\Psi_1=\Psi$
(see [MS]). Define 
$\v_t(x)=\frac{1}{a}\Psi_t(ax).$ Then $\v_t$ is a strict isotopy of $\Omega$
such that $\v_1(x)=\frac{1}{a} \Psi(ax)$.
For $t\in [0,1]$ write moreover
$\zeta_{t}(x)=\frac{1}{a(1-t)+t}\Psi((a(1-t)+t)x)$; then 
$\zeta_0=\v_1$ and $\zeta_1=\Psi$ and therefore
the composition of the isotopies
$\v_t$ and $\zeta_t$ define a strict isotopy of $\Omega$ as required.
{\bf q.e.d.}

\bigskip

Two proper
embeddings $\psi_1,\psi_2$ of 
a set $C\in {\cal O}$ into 
an open subset $\Omega$ of ${\bf R}^{2n}$ are
{\sl strictly isotopic} if there is a strict isotopy $\v_t$ of $\Omega$
such that $\v_1\psi_1\vert C=\psi_2\vert C$. As a corollary of Lemma 2.2
we obtain.

\bigskip
{\bf Corollary 2.3:} {\it Let $\Omega\subset \mathbb{R}^{2n}$ be 
open, bounded and star-shaped, and let $B\subset \mathbb{R}^{2n}$ be
open and bounded and such that 
$H^1(\overline{B},\mathbb{R})=0.$ 
Then any two proper embeddings
$\psi_1,\psi_2$
of $B$ into
$\Omega$ are strictly isotopic.}

{\it Proof:} The case $n=1$ is well known, 
so assume
that $n\geq 2$. Let $\Omega \subset \mathbb{R}^{2n}$ be open and
star-shaped with respect to the origin. Then $\overline{\Omega}$ is
simply connected and the same is true for
$\mathbb{R}^{2n}-\Omega$. Let 
$B\subset \mathbb{R}^{2n}$ be open and bounded and such that
$H^1(\overline{B},\mathbb{R})=0$.
Let
$\psi_1,\psi_2:B\to\Omega$ be proper
embeddings. Choose an
open neighborhood $U\supset {\overline B}$ of $B$
such that $\psi_i$ is
defined on $U$ $(i=1,2)$. Since $\mathbb{R}^{2n}-\Omega$ 
is simply connected and 
$H^1(\psi_i(\overline{B}),\mathbb{R})=0$ 
there is by
Theorem 2.1 a 
symplectomorphism $\Psi$ of $\mathbb{R}^{2n}$ whose
restriction to a neighborhood of $\mathbb{R}^{2n}-\Omega$ equals the
identity and whose restriction to a neighborhood of 
$\psi_1(\overline{B})$ 
which is contained in $\psi_1(U)$ 
coincides with $\psi_2\circ \psi_1^{-1}$.
Lemma 2.2 then shows that 
$\Psi$ is isotopic to the identity with an
isotopy which equals the identity on $\mathbb{R}^{2n}-\Omega$. 
{\bf q.e.d.}
\bigskip

{\bf Corollary 2.4:} {\it Let 
$\Omega\subset \mathbb{R}^{2n}$ be open, bounded and star-shaped and let
$U\subset V\in {\cal O}$ be such that
$\overline{U}\subset V$ and that
$\overline{U}$ and $\overline{V}$ are simply connected.
Let
$\v:U\to \Omega$ be a proper embedding. If there is 
a proper embedding $\zeta:V\to \Omega$ then
there is a proper embedding $\tilde \zeta:V\to \Omega$ whose
restriction to $U$ coincides with
$\v$.}

{\it Proof:} Let $\v:U\to \Omega$ and
$\zeta:V\to \Omega$ be proper embeddings. By Corollary 2.3
there is a symplectomorphism $\Psi$ of $\Omega$ which equals the 
identity near the boundary and such that
$\Psi\circ \v= \zeta\vert U$. Then $\Psi^{-1}\circ \zeta$ is
a proper embedding of $V$ whose restriction to
$U$ coincides with $\v$. 
{\bf q.e.d.}

\section{Cylindrical capacity and proper
displacement energy}

Using the assumptions and notations from the introduction and the
beginning of Section 2, the goal of this section is to show.

\bigskip
{\bf Theorem 3.1:} {\it There is
an open bounded starshaped set
$\Omega\subset \mathbb{R}^{4}$ such that
$c_p(\Omega)>e(\Omega).$}

\bigskip
For the proof of our theorem we 
will need the
following simple lemma.

\bigskip
{\bf Lemma 3.2:} {\it Let $h:\mathbb{R}^{2}\to \mathbb{R}$ and 
$f:\mathbb{R}^{2n-2}
\to \mathbb{R}$ be smooth functions with Hamiltonian flows
$\v_t,\eta_s$. View $h$ and $f$ as functions on $\mathbb{R}^{2n}$
which only depend on the first two and last $2n-2$ coordinates
respectively. Let $\nu_t$ be the Hamiltonian flow on 
$\mathbb{R}^{2n}$
of the function $h f$; then
$\nu_t(x,z)=(\v_{tf(z)}(x),\eta_{th(x)}(z))$ for every 
$x\in \mathbb{R}^2$
and every $z\in \mathbb{R}^{2n-2}$.

Proof:} Let $Z_h,Z_{f}$ be the Hamiltonian vector fields of $h,f$
as functions on $\mathbb{R}^{2n}$ only depending on the first two and
last $2n-2$ coordinates respectively. Then $Z_h$ is a section of
the $2$-dimensional subbundle of $T\mathbb{R}^{2n}$ spanned by the
basic vector fields $\frac{\partial}
{\partial x_1},\frac{\partial}
{\partial x_2}$, and $Z_f$ is a section of the $2n-2$-dimensional
subbundle of $T\mathbb{R}^{2n}$ spanned by the basic vector fields
$\frac{\partial}
{\partial x_i}$ for $i\geq 3$.
Moreover $f Z_h+hZ_{f}$ is the Hamiltonian vector field of the
function $h f$. We denote by $\nu_t$ its Hamiltonian flow.

The functions $h,f$ induce Hamiltonian flows $\v_t,\eta_t$ on
$\mathbb{R}^{2},\mathbb{R}^{2n-2}$. Since $h$ is constant along the
orbits of $\v_t$ and $f$ is constant along the orbits of $\eta_t$,
for every $x\in \mathbb{R}^2$ and every $z\in \mathbb{R}^{2n-2}$ we have
$\nu_t(x,z)=(\v_{tf(z)}(x),\eta_{th(x)}(z))$ which shows the lemma.
{\bf q.e.d.}
\bigskip

Using our lemma we can now determine the cylindrical
capacity of a special open bounded star-shaped set as 
follows. 

\bigskip

{\bf Example 3.3:} Consider $\mathbb{R}^4$ with the standard 
symplectic form $\omega_0$. Define 
$Q_1=\{(0,s,t,0)\mid -\frac{1}{2}\leq s\leq \frac{1}{2},0\leq t\leq 2\}$.
For a small number $\delta<1/4$ let $L$ be the convex cone
in the $(x_2,x_3)$-plane with vertex at the
origin whose boundary consists of the ray $\ell_1$ 
through $0$ and the point
$(0,\delta,2,0)$ and the ray $\ell_2$
through $0$ and $(0,\frac{1}{2},2,0)$. Let
$\tau >2$ be the unique number with the property that
the 
line 
$\{(0,s,\tau,0)\mid s\in 
\mathbb{R}\}$ intersects the cone $L$ in a segment of length $1$.
Define $Q_2=\{(0,x_2,x_3,0)\in L\mid
x_3\leq \tau\}$; then the boundary of
$Q_2$ is a triangle with one vertex at the origin, a 
second vertex
$z_1\not= 0$ on the line $\ell_1$ and the third vertex
$z_2\not=0$ on the line $\ell_2$.
Let $\ell_3$ be the line through
$z_2$ which is parallel to $\ell_1$.
The lines $\ell_1,\ell_3$ bound a strip $S$
which is foliated into line segments of length 1 which are
parallel to the $x_2$-coordinate axis.
Choose a large
number $M>\tau+2$ and define $Q_3=\{(0,x_2,x_3,0)\in S\mid
\tau\leq x_3\leq M\}$. 

By construction, the set $Q=Q_1\cup Q_2\cup Q_3$ is 
star-shaped with respect to $0$, and for every $t>0$
the line $\{x_3=t,x_1=x_4=0\}$ intersects
$Q$
in a connected segment of length at most $1$. 

Let $P_0\subset \mathbb{R}^3=\{x_1=0\}$ be the 
set which 
we 
obtain by rotating $Q$ about the origin in the
$(x_3,x_4)$-plane. The set
$P=[-1/2,1/2]\times P_0$ is star-shaped with
respect to the origin and it 
contains the cube $[-\frac{1}{2},
\frac{1}{2}]^2\times
D$ where $D$ is the disc of capacity $4\pi >1$ in the
$(x_3,x_4)$-plane. Thus the Gromov-width of $P$
is not smaller than $1$. 

We claim that the cylindrical capacity of $P$ equals $1$.
For this let $\epsilon >0$ and 
choose a smooth function
$\sigma:[0,\infty)\to [0,\infty)$
with the property that we have
$Q\subset \{(0,s,t,0)\mid t\geq 0,
\sigma(t)-1/2\leq s\leq \sigma(t)+1/2+\epsilon\}$.
Such a function $\sigma$
exists by the definition of $Q$, and 
we may assume that 
it vanishes
identically on $[0,2]$. 
The Hamiltonian flow $\zeta_t$ of the smooth function 
$(x_3,x_4)\to \sigma(\sqrt{x_3^2+x_4^2})$ preserves the concentric
circles about the origin. 
Define $h(x_1,x_2,x_3,x_4)=-x_1\sigma
(\sqrt{x_3^2+x_4^2}).$
By Lemma 3.2 and the fact that $P$ is invariant under
rotation about the origin in the $(x_3,x_4)$-plane we
conclude that
the image of $P$ under the time-one map
of the Hamiltonian flow of the function $h$ equals the set
$\hat P=\{(x_1,x_2,x_3,x_4)\mid (x_1,x_2+\sigma(\sqrt{x_3^2+x_4^2}),
x_3,x_4)\in P\}$ which is contained in the
subset $[-1/2,1/2+\epsilon]^2\times B^{2}(\pi M^2)$ 
of the cylinder
$[-1/2,1/2+\epsilon]^2\times \mathbb{R}^2\subset
\mathbb{R}^4$. 
Since $\epsilon >0$ was arbitrary,
the cylindrical capacity of 
$P$ is not bigger than $1$ and hence it coincides with
the Gromov width of $P$.

\bigskip

Now we can complete the
proof of Theorem 3.1.
Let $Q\subset \{x_1=0,x_4=0\}$
be as in Example 3.3. Reflect $Q$ along the line
$\{x_3=0\}$ in the $(x_2,x_3)$-plane. We obtain a set
$\tilde Q$ which is star-shaped with respect to
the origin. 
Define $\tilde P_0\subset \{x_1=0\}$ to be the set which 
be obtain by rotating $\tilde Q$ about the origin in
the $(x_3,x_4)$-plane.
Let
$\tilde P=[-1/2,1/2]\times \tilde P_0$. 
Then $\tilde P$ is star-shaped with respect to the origin
and contains $P$ as a proper subset.
 
We claim that $e(\tilde P)=1$. To see this notice that
for every $t> 0$ the intersection of $\tilde Q$
with the line 
$L_t=\{(0,s,t,0)\mid s\in \mathbb{R}\}$ 
consists of
at most $2$ segments of length at most $1$ each. 
Thus for every $\epsilon >0$ 
we can find a smooth function
$f_\epsilon$ on the 
half-plane $\{x_3>0\}$ in the 
plane $\{x_1=x_4=0\}$  
which satisfies 
$\sup_{z\in \tilde Q}f_\epsilon(z)-
\inf_{z\in \tilde Q}f_\epsilon(z)\leq 1+\epsilon$
and such 
that for
every $t\geq 0$ its restriction 
to each of the
at most two components of $L_t\cap
\tilde Q$ equals a 
translation.
We may choose $f_\epsilon$ in
such a way that $f_\epsilon(x_2,x_3)=x_2$ for $0<x_3<2$. 

Extend the function $f_\epsilon$ 
to a function
$f$ on $\mathbb{R}^4$
which does not depend on the first coordinate
and is invariant under rotation about the origin in the
$(x_3,x_4)$-plane. The Hamiltonian vector field of the
restriction of $-f$ to our set $\tilde P$ is of the
form $\frac{\partial}{\partial x_1}+Z$ where the vector
field $Z$ is tangent to the concentric circles 
about the origin in the
$(x_3,x_4)$-plane. Since 
$\tilde P$ is invariant under
rotation about the origin in the $(x_3,x_4)$-plane we
conclude that for every $s>0$ the image 
of $\tilde P$ under the
time-s map of the Hamiltonian flow of $ f$ equals
the set $[-1/2+s,1/2+s]\times \tilde P_0$. 
This means that 
the time-$(1+\epsilon)$ map
of the Hamiltonian flow of $f$ properly displaces
$\tilde P$. Via multiplying $f$ with a suitable
cutoff-function we deduce that 
the proper displacement energy of $\tilde P$ is not bigger
than $(1+\epsilon)^2$. 
Since $\epsilon >0$ was arbitrary and 
since $c_0(\tilde P)\geq 1$ we have
$e(\tilde P)=1$.

We are left with showing 
that
$c_p(\tilde P)>1$. 
For this define
$A\subset \mathbb{R}^2$ to be the closed annulus 
$\overline{B^2(\pi M^2)}-B^2(\pi \tau^2)$ of area 
$\pi(M^2-\tau^2)\geq  4\pi$.
By the discussion in Example 3.3 there is a 
small 
number $\rho>0$ depending on the choice of $\delta$
in the construction of the set $Q$ 
with the following properties.
\begin{enumerate}
\item
For every $\epsilon \in (0,\rho)$ 
there is a symplectic embedding $\psi_\epsilon$ of a neighborhood
of the  
star-shaped set
$P\subset \tilde P$ into the
cylinder
$B^2(1+\epsilon)\times \mathbb{R}^2$ with the property
that $\psi_\epsilon(P)\supset
B^4(1-\epsilon)\cup B^2(1-\epsilon)\times A$ and
$\psi_\epsilon(P\cap \{x_2>0\})\supset
B^2(\rho)\times B^2(\pi M^2)$.
\item
There is a proper symplectic embedding of
the standard ball
$B^4(1/2-\rho)$ into 
$\tilde P-P$  
whose image $\tilde B$
is strictly isotopic in $\tilde P\cap \{x_2<0\}\supset
\tilde P-P$
to a standard ball embedded in $B^4(1-\epsilon)\cap \{x_2<0\}$.
\end{enumerate}

Assume to the
contrary that $c_p(\tilde P)=1$. 
Then there is for every
$\epsilon \in (0,\rho)$ a proper symplectic embedding 
of $\tilde P$
into the cylinder 
$B^2(1+\epsilon)\times \mathbb{R}^2$. 
Since the closure of the disconnected set $P\cup \tilde B\subset
\tilde P$ is simply connected and since the
standard cylinder $B^2(1+\epsilon)\times \mathbb{R}^2$ 
is star-shaped with respect to the origin, we
can apply 
Corollary 2.4 to proper embeddings of 
$P\cup \tilde B$
into $B^2(1+\epsilon)\times \mathbb{R}^2$.
This means that there is a proper
symplectic embedding $\Psi$ of $\tilde P$ into
$B^2(1+\epsilon)\times \mathbb{R}^2$ whose
restriction to $P\subset \tilde P$ coincides
with $\psi_\epsilon$ and which
maps $\tilde B$ to a standard euclidean ball 
$\hat B$ which is 
contained in
$B^2(1+\epsilon)\times \bigl(\mathbb{R}^2-B^2(\pi M^2)\bigr)$
and can be obtained from $B^4(1/2-\rho)$ by a translation.

Let 
$\omega_1$ be a standard volume form
on the sphere $S^2$ whose total area is bigger than but
arbitrarily close to 
$1+\epsilon$.
Embed the disc $B^2(1+\epsilon)$
symplecticly into $S^2$.
The image in $S^2$ of the annulus
$B^2(1+\epsilon)-B^2(1-\epsilon)$ is contained in a closed round
disc $D\subset S^2$ of area bigger than but arbitrarily close
to $2\epsilon$. The complement of $D$ in $S^2$ is the disc
$B^2(1-\epsilon)$. The complement of the disc $B^2(\rho)$
in $S^2$ is area-preserving equivalent to a 
closed disc in $\mathbb{R}^2$.

Our embedding of $B^2(1+\epsilon)$ into $(S^2,\omega_1)$
extends to a symplectic embedding of $B^2(1+\epsilon)\times 
\mathbb{R}^2$ into $(S^2\times \mathbb{R}^2,
\omega=\omega_1+\omega_0)$.
Thus if $\tilde P$ admits a 
proper symplectic embedding
into $B^2(1+\epsilon)\times \mathbb{R}^2$ then 
the standard linear embedding of $B^4(1/2-\rho)$ onto a ball in
$S^2\times
(\mathbb{R}^2-B^2(\pi M^2))\subset
S^2\times \mathbb{R}^2$ is strictly isotopic in
$S^2\times (\mathbb{R}^2-B^2(\pi M^2))
\cup D\times A\cup B^2(1+\epsilon)\times 
B^2(\pi \tau^2)\subset S^2\times \mathbb{R}^2$ 
to the standard
inclusion of $B^4(1/2-\rho)$
into $B^2(1+\epsilon)\times B^2(1/2-\rho)$.

However 
a suitable version of the 
{\sl symplectic camel
theorem} in dimension 4 
[MDT] shows that
for sufficiently small $\epsilon$ this is not possible. 
We formulate this version as a
proposition which 
completes the proof of our Theorem 3.1.

\bigskip

{\bf Proposition 3.4:} {\it  Let 
$D\subset S^2$ be an open round
disc of area $\frac{1}{4}$ in a standard sphere 
$(S^2,\omega_1)$ of area $1$. 
Let $\omega=\omega_1
+\omega_2$ be a standard symplectic form on $S^2\times 
\mathbb{R}^2$.
For every $R\in (\frac{1}{4},1)$ a standard 
embedding of the
ball $B^4(R)$ into 
$S^2\times (\mathbb{R}^2-B^2(3))$
is not properly isotopic in
$S^2\times (\mathbb{R}^2-\overline{B^2(2)})\cup
D\times \partial{B^2(2)}\cup
B^2(1)\times B^2(2)\subset S^2\times \mathbb{R}^2$
to a standard embedding of $B^4(R)$
into
$S^2\times B^2(1)$.}

{\it Proof:} Using the notation from the proposition,
let $R\in (\frac{1}{4},1)$ and 
let $\v_0$ be a standard embedding of the ball $B^4(R)$ of
capacity $R$ into $B^2(R)\times 
(\mathbb{R}^2-B^2(3))\subset S^2\times (\mathbb{R}^2-B^2(3))$.
Let moreover $\v_1$ be a standard embedding of $B^4(R)$ into
$B^2(1)\times B^2(1)\subset S^2\times B^2(1)$. 
We argue by contradiction and we assume that 
$\v_0$ can be connected to $\v_1$ by a proper isotopy 
$\v_t$ $(t\in [0,1])$ 
whose image is contained in
the subset 
$S^2\times (\mathbb{R}^2-\overline{B^2(2)})\cup
D\times \partial B^2(2)\cup B^2(1)\times B^2(2)$
of the manifold $S^2\times \mathbb{R}^2$.

We follow
[MDT] and arrive at a contradiction
in three steps.

\smallskip

{\sl Step 1}

Let $c$ be the boundary of the disc $B^2(2)$.
Denote by $\nu$ the boundary circle of 
the disc $D\subset S^2$.
Then $T=\nu\times c$
is a Lagrangian torus embedded   
in $S^2 \times \mathbb{R}^2$. 
For a fixed point $y$ on $c$,
the circle
$\nu\times \{y\}$
bounds the embedded disc 
$D\times \{y\}\subset S^2\times \{y\}$.
We call such a disc a {\sl standard flat disc}.
It defines
a homotopy class of maps of pairs from
a closed unit disc
$(D_0,\partial D_0)\subset \mathbb{R}^2$ 
into $(S^2\times \mathbb{R}^2,T)$.

Let ${\cal J}$ be the space of all smooth 
almost complex structures $J$ 
on $S^2\times \mathbb{R}^2$
which {\sl calibrate} 
the symplectic form $\omega$ 
(i.e. such that $g(v,w)=\omega(v,Jw)$ defines
a Riemannian metric on $S^2\times \mathbb{R}^2$).
In the sequel we mean by a pseudoholomorphic disc 
a disc which is holomorphic with respect to some
structure $J\in {\cal J}$.
For $J\in {\cal J}$ 
define a {\sl $J$-filling} of the torus $T$ to be
a $1$-parameter family of disjoint, $J$-holomorphic
discs 
which are homotopic 
as maps of pairs 
to the standard flat disc, 
whose boundaries foliate $T$ and whose union $F(J)$
is diffeomorphic to $D\times S^1$ and does not
intersect $(S^2-D)\times c$.
The set $F(J)$ then necessarily disconnects 
$\Omega=S^2\times \mathbb{R}^2-(S^2-D)\times c$.

For $t\in [0,1]$ let $J_t\in {\cal J}$ be an almost
complex structure depending continuously on $t$. We require that
$J_0=J_1$ is the standard complex 
structure and that the restriction
of $J_t$ to $\v_tB^4(R)$ coincides with 
$(\v_t)_* J_0$ (where by
abuse of notation we denote by $J_0$ the natural
complex structure on $\mathbb{R}^4$ and
on $S^2\times \mathbb{R}^2$). Such structures exist since
${\cal J}$ is the space of smooth sections of a fibre
bundle over $S^2\times \mathbb{R}^2$ with
contractible fibre.

Assume that for every $t\in [0,1]$ there
is a unique $J_t$-filling $F(J_t)$ of $T$ depending
continuously on $t$ in the Hausdorff topology for
closed subsets of $S^2\times \mathbb{R}^2$. 
Then 
\[X=\{(t,x)\mid x\in F(J_t),0\leq t\leq 1\}\]
is a closed subset of $[0,1]\times \Omega$. 
Since each filling $F(J_t)$ 
disconnects $\Omega$, the set $X$
disconnects $[0,1]\times \Omega$. 
Now $J_0$ and $J_1$ are standard
and the standard filling 
of $T$ by flat discs 
separates 
$\v_0B^4(R)$ from $\v_1 B^4(R)$.
Thus the points
$(0,\v_0(0))$ and $(1,\v_1(0))$ are contained in different components
of $[0,1]\times \Omega-X$. 
Therefore the path $(t,\v_t(0))$ must intersect
$X$. In other words, for some $t$, there is a $J_t$-holomorphic
disc $C$ through $\v_t(0)$
with boundary on $T$ 
and which is contained 
in the homotopy class of the standard
flat disc.
The connected component of 
$\v_t^{-1}C\cap  
B^4(R)$ containing $0$ 
is a planar holomorphic curve with respect
to the standard integrable complex structure
whose boundary is contained
in the boundary of $B^4(R)$. 
Thus this surface is 
a minimal surface with boundary on the boundary
of $B^4(R)$ and therefore its area
is not smaller than $R$ [G]. On the other hand,
the torus $T$ is Lagrangian and hence the
area of a $J_t$-holomorphic disc with
boundary on $T$ only depends on the free
homotopy class of the boundary curve. In particular,
since the curve $C$ is 
homotopic to the standard
flat disc its area
equals the area $\frac{1}{4}$ of the standard flat disc.
This contradicts our assumption that
$\frac{1}{4}<R$
(compare [MDT] p.178). 
 
By the above it is now enough to construct
for every 
$t\in [0,1]$ an almost complex structure $J_t\in
{\cal J}$ whose restriction to $\v_tB^4(R)$ coincides
with $(\v_t)_*J_0$ and such that 
for every $t\in [0,1]$ 
the torus $T$
admits a unique $J_t$-filling depending
continuously on $t\in [0,1]$ in the Hausdorff topology. 
For this we follow again [MDT]. 

\smallskip

{\sl Step 2}

Let $H$ be an oriented hypersurface in
an almost complex $4$-manifold $(N,J)$.
There is a unique two-dimensional subbundle
$\xi$ of the tangent bundle of $H$ which is
invariant under $J$.  
We call $H$
{\sl $J$-convex} if for one (and hence any)
one-form $\alpha$ on $H$ whose kernel 
equals 
$\xi$ and which defines together with the restriction
of $J$ to $\xi $ the orientation of $H$  
and for every
$0\not=v\in \xi$ we have $d\alpha(v,Jv)> 0$.

Let $H\subset S^2\times
\mathbb{R}^2$ be a smooth hypersurface which contains
the Lagrangian torus $T$ and bounds
an open domain $U_H\subset
S^2\times \mathbb{R}^2$ 
which contains the
open disc bundle $D\times c-T$. We equip
$H$ with the orientation induced by the outer normal
of $U_H$. 
Denote by ${\cal J}_H\subset {\cal J}$
the set of all almost complex structures $J\in 
{\cal J}$ for which $H$ is $J$-convex and which
coincide with the standard complex structure $J_0$
near $T$. If the set ${\cal J}_H$ is not empty
then it follows again from the fact that ${\cal J}$ is
the space of smooth sections of a fibre bundle over
$S^2\times \mathbb{R}^2$ with contractible fibre that
we can find some
$\tilde J\in {\cal J}_H$ 
which coincides with $J_0$ on a neighborhood of the disc bundle
$D\times c$. Then the standard flat discs define a 
$\tilde J$-filling of $T$.

Let $J_t\subset {\cal J}_H$ $(t\in [1,2])$ be a differentiable
curve (there is some subtlety here about the differentiable
structure of ${\cal J}_H$ which will be ignored in
the sequel; a discussion of this
problem is contained in [MDT]). 
Assume that for every $t\in [1,2]$ there is a holomorphic
disc $D_t$ of $J_t$ with boundary on $T$ depending
continuously on $t$ 
in the Hausdorff topology 
and such that $D_1$ is a standard flat
disc. 
Since $H$ is $J_t$-convex,
Lemma 2.4 of [MD2] and 
Proposition 3.2 in [MDT] show that
each of the discs $D_t$ meets 
the hypersurface
$H$ transversely at its boundary, and its
interior
is contained
in $U_H$.

By assumption, each of the structures $J_t$ coincides
with the standard complex structure near the torus $T$. Thus there
is an open neighborhood $V$ of $T$ 
in $S^2\times \mathbb{R}^2$ 
and a $J_t$- antiholomorphic
involution in $T$ on $V$ (see [MDT]). 
Namely, under the usual identification of $\mathbb{R}^4$ with
$\mathbb{C}^2$, the torus $T$ is just the 
cartesian product of the boundary of a disc of radius $r_1>0$
with the boundary of a disc of radius $r_2>0$. 
The map $(z_1,z_2)\to (\frac{r_1^2}{\overline {z}_1},
\frac{r_2^2}{\overline{z}_2})$ is an antiholomorphic reflection
in $\mathbb{C}^2$ which fixes $T$ pointwise.
This reflection restricts to a 
reflection on an open neighborhood $V$ of $T$ which
is antiholomorphic with respect to any of the almost
complex structures $J_t$. 

Using this involution 
we can double
our domain $U_H$ near $T$ [MDT] and use 
the Schwarz reflection principle to extend our holomorphic
discs $D_t$ to holomorphic spheres in the double of $U_H$ [MDT]. 
Intersection theory
for pseudoholomorphic spheres in almost complex manifolds
then shows that each of our discs $D_t$ is embedded. Moreover
any two different
such discs for the same structure $J\in {\cal J}_H$ 
do not intersect [MD1].

Now Gromov's compactness theorem is valid for 
pseudoholomorphic discs with Lagrangian boundary 
condition [O]. The area of each 
pseudoholomorphic disc with boundary on $T$ and in 
the homotopy class of the standard flat disc
coincides with the area $\frac{1}{4}$ of the standard flat
disc. Moreover, since our torus $T$ is rational [P] and
$\frac{1}{4}$ is the generator of the subgroup of $\mathbb{R}$
induced by evaluation of $\omega_0$ on $\pi_2(\mathbb{R}^4,T)$,
$\frac{1}{4}$ is the minimal area of any 
pseudoholomorphic disc
whose boundary is contained in $T$.
This implies that for discs in our given homotopy class,
bubbling off of holomorphic spheres and
holomorphic discs can not occur. Therefore we can use
standard Fredholm theory for the Cauchy Riemann operator 
[MDT]
to compute for a dense set of points in
${\cal J}_H$ the parameter space of pseudoholomorphic discs
with boundary on $T$ and which are homotopic to the standard
flat disc. As a consequence [MDT],
for every $J\in {\cal J}_H$ which can be connected to
our fixed structure $\tilde J$ 
by a differentiable curve in ${\cal J}_H$ there
is 
a unique $J$-filling $F(J)$ of $T$ which depends continuously
on $J\in {\cal J}_H$ 
in the Hausdorff topology for closed subsets of $S^2\times
\mathbb{R}^2$ (here uniqueness means uniqueness of the
image and hence we divide the family of all holomorphic discs
by the group of 
biholomorphic automorphisms
of the unit disc in $\mathbb{C}$).

Together with Step 1 above we conclude that
our proposition follows if we can construct
a hypersurface $H$
in $S^2\times \mathbb{R}^2$ with the following properties.
\begin{enumerate}
\item $H$ contains $T$ and
bounds an open 
set $U_H$ containing 
the open disc bundle $D\times c-T$.
\item  $U_H$ contains a neighborhood 
of $\cup_{t\in [0,1]}\v_t B^4(R).$
\item 
${\cal J}_H\not=\emptyset.$
\end{enumerate}

Namely, for such a hypersurface $H$ we can choose
a fixed almost complex structure $J\in {\cal J}_H$ whose
restriction to the disc bundle $D\times c$ coincides with the
standard structure. For each $t\in [0,1]$ we modify 
$J$ near $\v_tB^4(R)$ in such a way that the modified
structure $J_t$ 
is unchanged near
the hypersurface $H$
and coincides with
$(\v_t)_*J_0$ on
$\v_tB^4(R)$.
We can do this in such a way that
$J_t$ depends differentiably on $t$. 
For each of the structures
$J_t$ there is then a unique $J_t$-filling of $T$ depending
continuously on $t$. 

By our assumption, the closure 
of the set $\cup_{t\in [0,1]}\v_t
B^4(R)$ is contained in $S^2\times 
(\mathbb{R}^2-\overline{B^2(2)})\cup D\times
\partial B^2(2)\cup B^2(1)\times B^2(2)$ and
hence 
it is enough to find a hypersurface $H$
in the symplectic manifold 
\[N=S^2\times (\mathbb{R}^2-\overline{B^2(2)})\cup
B^2(1)\times\partial B^2(2)\cup
\mathbb{R}^2\times B^2(2)\] with properties 1-3.
In the third step
of our proof we construct such a hypersurface.

\smallskip

{\sl Step 3}

Let $A\subset \mathbb{R}^2$ be a closed circular annulus
containing the circle $c$ in its interior. 
We assume that $A$ is small 
enough so that the closure of the set
$\cup_{t\in [0,1]}\v_tB^4(R)$ intersects
$S^2\times A$ in $B^2(\frac{1}{4})\times A$.
We also require that there is
some $a>0$ such that
$S^2\times A\subset S^2\times \mathbb{R}^2$ 
is symplectomorphic to the quotient 
of the
bundle $S^2\times [-a,a]\times \mathbb{R}$ 
under a translation $\tau$ in the plane $\mathbb{R}^2$ 
in such a way that the
torus $T$ is the quotient of the standard circle
bundle $\partial B^2(\frac{1}{4})\times \{0\}\times \mathbb{R}$.
The standard complex structure 
on $S^2\times \mathbb{R}^2$
is invariant under the translation $\tau$ and projects
to an integrable complex structure 
on a neighborhood of 
$S^2\times A$ in
$S^2\times \mathbb{R}^2$ which 
calibrates $\omega$.

For small $\sigma \in [0,a)$ 
write 
$\ell_\sigma=
\{(0,0,-\sigma,s)\mid s\in \mathbb{R}\}\subset \mathbb{R}^4$.
The circle
bundle $\partial B^2(\frac{1}{4})\times \{0\}\times \mathbb{R}$ 
is contained in the boundary 
$\partial U_\sigma$ 
of a tubular neighborhood $U_\sigma$ 
of some radius $r(\sigma)>1/2\sqrt{\pi}$ 
about the line $\ell_\sigma$. 
The hypersurface 
$\partial U_\sigma$ with
its orientation as the boundary of $U_\sigma$ 
is convex with respect to the euclidean 
metric. 

Let $X_\sigma$ be the 
gradient of the function $z\to \frac{1}{2}\,
{\rm dist}(z,\ell_\sigma)^2.$
Then $X_\sigma$ is perpendicular to the hypersurface
$\partial U_\sigma$ and can be
written down explicitly by
$X_\sigma= x_1\frac{\partial}{\partial x_1}
+x_2\frac{\partial}{\partial x_2}+(x_3+\sigma)\frac{\partial}{\partial x_3}$.
If we denote
by $\iota_{X_\sigma}\omega_0$ 
the 1-form $\omega_0(X_\sigma,\cdot)$ then
$\iota_{X_\sigma}=x_1dx_2-x_2dx_1+(x_3+\sigma)dx_4$ and
$d(\iota_{X_\sigma}\omega_0)=2dx_1\wedge dx_2+dx_3\wedge dx_4$.
This implies that
$\partial U_\sigma$ 
is 
convex with respect to 
$J_0$. 

As a first step towards the construction of a
hypersurface $H$ with properties 1-3 above
we construct a smooth embedded hypersurface
$E\subset \mathbb{R}^4-\ell_\sigma$  
with the following properties.
\begin{enumerate}
\item[a)] $E$
divides $\mathbb{R}^4$ into two components and is invariant
under the translations in direction of the $x_4$-axis.
\item[b)] There is some $\sigma\in (0,a/2)$ such that $E$
is everywhere transverse to 
the vector field $X_\sigma$ 
and contains a neighborhood of $T$ in $\partial U_\sigma$.
\end{enumerate}
For this let $E_0$ be a
hypersurface
which is contained in $\{0\leq x_3\leq a/2\}$ 
and 
which bounds a noncompact convex set $V$
containing $\{x_3>a/2\}$. We assume that $E_0$ is
invariant under translations 
along the lines parallel to $\ell_\sigma$, and
we orient $E_0$ as the boundary of $V$.

Assume that the
intersection of $E_0$ with the hyperplane
$\{x_3=0\}$ equals the line $\ell_0=\{(0,0,0,s)\mid s\in
\mathbb{R}\}$ and that
the hypersurface $E_0$ is
strictly convex in directions transverse
to the lines parallel to the $x_4$-axis. 
Since the vector field $X_0$ is tangent to the lines in
$\mathbb{R}^4$ which meet $\ell_0$
orthogonally, by convexity 
$X_0$ is 
everywhere transverse to $E_0-\ell_0$. 
More precisely, 
for every $z\in E_0-\ell_0$ the 
nonoriented angle between $X_0(z)$ and
the outer normal of $E_0$ (as the boundary of $V$)
at $z$ is smaller than $\pi/2$.
An explicit example of such a hypersurface can be
obtained as follows. Choose an even strictly convex function
$f:\mathbb{R}\to [0,a/2]$ which has
a (necessarily unique) minimum $0$ at $0$ and define
$E_0$ to be the solution of the equation
$f(x_1^2+x_2^2+x_4^2)-x_3=0.$

By construction, there is a tubular neighborhood 
of the line $\ell_0$ which is invariant under the
translation $\tau$ and contained in each
of the sets $U_\sigma$ for all small $\sigma \geq 0$. 
This implies that for a suitable choice of $E_0$ and 
for sufficiently
small $\sigma >0$ the vector field $X_\sigma$ is
everywhere transverse to $E_0-U_\sigma$. 
More precisely, for such a $\sigma$  
and for every compact set $K$ there is
a number $\delta(K)>0$ such that for
$z\in K\cap (E_0-U_\sigma)$ the nonoriented angle at $z$ between
$X_\sigma$ and the \emph{outer} normal of
$E_0$ is contained in the interval $[0,\pi/2-\delta(K)]$.
This means that for such a $\sigma$ 
the
vector field $X_\sigma$ is everywhere transverse to
the boundary of the set
$U_\sigma\cup V\cup \{x_3<-\sigma/2\}=W$, and its nonoriented
angle with the outer normal
of this boundary is stricly smaller than
$\pi/2$. In particular, after a small
perturbation of our sets near the intersections 
$\partial U_\sigma\cap E_0$ and $\partial U_\sigma\cap
\{x_3=-\sigma/2\}$ 
we may assume that the boundary of
our set $W$ is a smooth hypersurface $\partial W$ which is
contained in $\{-\sigma/2\leq x_3\leq a/2\}$,
is invariant under
the translation $\tau$ and 
satisfies properties a),b) above.

The
boundary $\partial W$ of $W$ is contained in 
$\{-\sigma/2\leq x_3\leq
a/2\}$ and hence it 
projects to a smooth hypersurface 
in $\mathbb{R}^2\times A$ 
which bounds an open set $\tilde W$.
The set $\tilde W$ in turn projects to a set $\hat W\subset N$.
By our explicit construction we may assume that
$\hat W$ 
contains the closure of the set $\cup_{t\in [0,1]}
\v_tB^4(R)$.

The vector field $X_\sigma$ is transverse
to $\partial W$ and therefore 
the kernel of the 
$1$-form $\iota_{X_\sigma}\omega_0$ intersects
the tangent
bundle of $E$ 
in a two-dimensional subbundle $\xi$.
The restriction of $\omega_0$ to $\xi$
is non-degenerate and hence 
there is an almost complex structure 
$\tilde J_\sigma$ on
$\xi$ which calibrates $\omega_0\vert \xi$. 
This structure $\tilde J_\sigma$ extends to an almost
complex structure $J_\sigma$ near $E$ which calibrates
$\omega_0$. 
Since
a neighborhood of $T$ in $E$ is contained in $\partial U_\sigma$,
the almost complex structure $J_\sigma$ 
can be chosen to coincide with
the standard complex structure near $T$. 
Moreover, 
$E$ is $J_\sigma$-convex.
Now $X_\sigma$ and $\partial W$ are
invariant under the translation $\tau$ 
and hence we may assume that the same is true for the 
almost complex structure $J_\sigma$.
Therefore this 
almost complex structure 
projects to an almost complex structure 
on a neighborhood of $\partial \tilde W$ 
which we denote again by $J_\sigma$. 
The hypersurface $\partial\tilde W$ is $J_\sigma$-convex.

By construction, 
there is a circle $\gamma$ in 
$\mathbb{R}^2-\overline{B^2(2)}$ such that
the intersection of $\partial \tilde W$ with the
set $(\mathbb{R}^2-B^2(1/2))\times (\mathbb{R}^2-B^2(2))$
is contained in the  hypersurface 
$\mathbb{R}^2\times \gamma$.
But this just means
that $\partial\tilde W$ projects to a smooth hypersurface $H$ in
the set $N$.

Denote by
$\alpha$ the restriction of the
1-form $\iota_{X_\sigma}\omega_0$ to 
the hyperplane
$Q=\{x_3=-\sigma/2\}$.
Let $(r,\theta)$ be polar coordinates
about $0$ in the $(x_1,x_2)$-plane.
By our explicit formula for $X_\sigma$ the form
$\alpha$ 
can be written in coordinates
$(r,\theta,x_4)$ for $Q$ as
$\alpha=rd\theta +\frac{\sigma}{2}dx_4$. Then
$d\alpha=dr\wedge d\theta$.

Now if $\beta$ is any one-form
on $\{x_3=-\sigma/2\}$ which 
is invariant under rotation 
about the origin 
in
the $(x_1,x_2)$-plane and the
translations along the lines parallel to $\ell_\sigma$ 
then in our above coordinates we can write 
$ \beta=\varphi(r)d\theta+ \rho(r)dx_4$ for 
functions $\varphi,\rho$ on $(0,\infty)$. The
one-form $\beta$ vanishes nowhere if and only
if the functions $\varphi,\rho$ do not have
a common zero.
Now $d\beta=
\varphi^\prime(r)dr\wedge d\theta +\rho^\prime(r)dr\wedge dx_4$
and hence
the restriction of $d\beta$ to the kernel of $\beta$ 
vanishes nowhere if and only if
we have
$\varphi^\prime(r)\rho(r)-\rho^\prime(r)\varphi(r)>0$.

Let $\varphi, \rho$ be functions
on $(0,\infty)$ 
with the following properties.
\begin{enumerate}
\item $\varphi(t)=t, \rho(t)=\sigma/2$ for $t\leq 1/2$.
\item $\varphi^\prime  \rho-\rho^\prime \varphi>0.$
\item $\varphi$ and $\rho$ do not have a common zero.
\item $\varphi(1)=0.$
\end{enumerate}
Such functions $\varphi,\rho$ can easily be constructed.
We replace the one-form $\alpha$ on $Q$ by the one-form
$\beta=\varphi(r)d\theta
+\rho(r)dx_4$. Then $\beta$ is a contact form
on $Q$ which coincides with the contact
form $\alpha$ on $B^2(1/2)\times \mathbb{R}^2\cap Q$. Moreover,
this contact form projects to a contact form on the
hypersurface $\mathbb{R}^2\times \gamma$ which we
denote again by $\beta$.
The kernel of $\beta$   
on $\partial B^2(1)\times \gamma$ 
is spanned by
the vector fields $\frac{\partial}{\partial x_1},
\frac{\partial}{\partial x_2}$ and consequently 
this kernel projects to a 
plane bundle on the projection of 
$(\overline{B^2(1)}-B^2(1/2))\times \gamma$ to $N$ which
we obtain by mapping each circle $\partial B^2(1)\times \{y\}$
to a point. 

In other words, there is a modification 
$\tilde J_\sigma$
of the almost complex
structure $J_\sigma$ which coincides with the integrable
complex structure near the torus $T$ 
and which projects  
to an almost complex structure on 
a neighborhood of $H$ in $N$.
This structure is the restriction of  
a smooth almost complex structure $\hat J_\sigma$
on $N$ which calibrates $\omega$.
The oriented hypersurface $H$ is the boundary
of an open set $U$ containing
$\cup_{t\in [0,1]}\v_t B^4(R)$, and it
is $\hat J_\sigma$-convex.

Together this means that $H$ and $\hat J_\sigma$ satisfy the
properties 1-3 above. This finishes the proof of our proposition.
{\bf q.e.d.}

\section{Displacement and proper displacement}
\bigskip

In this short section we collect some easy
properties of the displacement energy and
the proper displacement energy.
We continue to use the
assumptions and notations from Section 2-4.

Recall that the {\sl displacement energy} $d(\Omega)$
of an open bounded set $\Omega\subset \mathbb{R}^{2n}$ 
equals
the
infimum of the Hofer norms of all symplectomorphisms $\Psi$ of
$\mathbb{R}^{2n}$ such that $\Psi\Omega\cap \Omega= \emptyset$. 
For every $\Omega\in {\cal O}$ the displacement
energy $d(\overline{\Omega})$ of the closure
$\overline{\Omega}$ of $\Omega$ 
is not smaller
than the proper displacement energy $e(\Omega)$ of $\Omega$.
In the case $n=1$ equality $d(\overline{\Omega})=e(\Omega)$ 
always holds.

The following lemma is an
easy corollary of the neighborhood extension
theorem of Banyaga [B] (see Theorem 2.1).

\bigskip

{\bf Lemma 4.1:} {\it Let  $\Omega\in {\cal O}$
be such that $H^1(\overline{\Omega},\mathbb{R})=0$; then
$d(\overline{\Omega})=e(\Omega)$.}

{\it Proof:} Since we always have
$e(\Omega)\geq d(\overline{\Omega})$ we
have to show the reverse inequality under the 
assumption that $H^1(\overline{\Omega},\mathbb{R})=0.$
For this
we only have to consider the case $n\geq 2$.

Let $\epsilon >0$ and
let $\Psi\in {\cal D}$ be a
compactly supported Hamiltonian
symplectomorphism of $\mathbb{R}^{2n}$
of Hofer norm smaller than $d(\overline{\Omega})+\epsilon$
and such that $\Psi(\overline{\Omega})\cap
\overline{\Omega}=\emptyset$. Then there is an open
neighborhood $U$ of $\overline{\Omega}$ such that
$\Psi(U)\cap U=\emptyset.$

Assume without loss of generality that $U\subset \{x_1<0\}$. Let
$e_1$ be the first basis vector of the standard basis of 
$\mathbb{R}^{2n}$, choose some 
$\mu<\inf\{x_1(z)\mid z\in U\}$ and define
$W=U\cup (U-\mu e_1).$ Then $W$ contains two copies of
$\overline{\Omega}$ in its interior which are separated by the
hyperplane $\{x_1=0\}$. Moreover the set $\Omega\cup \Psi\Omega$
admits a natural proper symplectic embedding into $W$ whose
restriction to $\Omega$ is just the inclusion.

Since $H^1(\overline{\Omega}\cup \overline{\Psi\Omega},
\mathbb{R})=0$,
by the Banyaga extension theorem this proper symplectic
embedding of $\Omega\cup \Psi \Omega$ into $W$ can be extended
to
a symplectomorphism $\eta$ of $\mathbb{R}^{2n}$.
Then $\eta\circ
\Psi\circ \eta^{-1}$ is a symplectomorphism of Hofer-norm smaller
than $d(\overline{\Omega})+\epsilon$ which properly displaces
$\Omega=\eta(\Omega)$. This shows that $e(\Omega)\leq
d(\overline{\Omega})+\epsilon$, and since $\epsilon >0$ was
arbitrary the lemma follows. {\bf q.e.d.}

\bigskip

Finally we look at
the relation
between  the displacement energy of a set $\Omega\in {\cal
O}$ and the displacement energy of its closure
$\overline{\Omega}$. We first give an easy example which shows
that the equality $d(\Omega)=d(\overline{\Omega})$ 
does not even hold for
open bounded topological balls with smooth boundary.

\bigskip

{\bf Example 4.2:} Let $\epsilon \in (0,1/2)$
and define \[Q_\epsilon =
\overline{B^2(4)}\times [0,1]^2\cup
(\overline{B^2(8)}-B^2(4))\times
[1-\epsilon,2-\epsilon]\times [0,1].\] Then $Q_\epsilon$ is a
closed topological ball with piecewise smooth boundary whose
interior we denote by $U_\epsilon$. By construction, the sets
$U_\epsilon$ and $U_\epsilon +(0,0,1,0)$ are disjoint and
hence $U_\epsilon$ can be displaced by the
time-one map of the Hamiltonian flow of the function
$f(x_1,x_2,x_3,x_4)=-x_4$. This implies that $d(U_\epsilon)\leq 1$.
Since $U_\epsilon$ contains the
open cylinder $B^2(4)\times
(0,1)^2$ of displacement energy $1$ we conclude that 
$d(U_\epsilon)=1$.

On the other hand, the closure $Q_\epsilon$ of $U_\epsilon$
contains the split Lagrangian torus $T^2=\partial B^2(4)\times
\partial([0,2-\epsilon]\times [0,1])$ of displacement energy
$2-\epsilon$ and therefore we have $d(Q_\epsilon)>d(U_\epsilon)$.
Via replacing the squares in our construction by discs with smooth
boundary we can also find an example of a topological ball
$\Omega$ with smooth boundary $\partial \Omega$ and such that
$d(\overline{\Omega})>d(\Omega)$.

\bigskip

Recall that a (not necessarily smooth) hypersurface $H$ in 
$\mathbb{R}^{2n}$
is of {\sl contact type} if there is a conformal vector
field $\xi$ defined near the hypersurface $H$ (i.e. such that the
Lie-derivative of $\omega_0$ with respect to $\xi$ coincides with
$\omega_0$) which is transverse to $H$ in the sense that the flow
lines of $\xi$ intersect $H$ transversely. Let $\v_t$ be the
local flow of $\xi$ and assume that there is some $\epsilon >0$
such that $\v_t$ is defined near $H$ for all $t\in (-\epsilon,
\epsilon)$.
If $H$ is
the boundary of a bounded open set $\Omega\in {\cal O}$
then for $t\in (0,\epsilon)$ the set
$\v_tH$ is the boundary of a neighborhood of $\overline{\Omega}$,
and for $t\in (-\epsilon,0)$ the set $\v_t H$ is contained in
$\Omega$. The hypersurface $H$ is called of {\sl restricted
contact type} if the conformal vector field $\xi$ can be defined
on all of $\mathbb{R}^{2n}$. For example, if $\Omega\in {\cal O}$ is
starshaped with respect to $0$ and if the lines 
through $0$ intersect the boundary $\partial \Omega$
transversely then $\partial \Omega$ is of restricted contact type.

Our last lemma shows that the difficulty encountered in our 
example 4.2 does not occur for open sets with boundary of restricted
contact type. 

\bigskip
{\bf Lemma 4.3:} {\it Let $\Omega$ be an open bounded set
in $\mathbb{R}^{2n}$ $(n\geq 2)$. If the boundary of $\Omega$ is
of restricted contact type then $d(\Omega)=d(\overline{\Omega})$.}

{\it Proof:} Let $\Omega\in {\cal O}$ be an open bounded subset of
$\mathbb{R}^{2n}$ whose boundary is of
restricted contact type.

Let $\xi$ be a conformal vector field on $\mathbb{R}^{2n}$ which
intersects the boundary of $\Omega$ transversely.
Assume that there is a neighborhood $U$ of $\Omega$ and a number
$\epsilon >0$ such that the local
flow $\v_t$ of $\xi$ is defined on $(-2\epsilon,2\epsilon)\times U$.
Then
the image of $\Omega$ under the time-$\epsilon$ map
of the flow $\v_t$ of $\xi$ is a neighborhood of
$\overline{\Omega}$. Since $\v_t^*\omega_0=e^t\omega_0$ for all
$t$ and wherever this is defined, by conformality the displacement
energy of $\v_\epsilon \Omega$ is not bigger than $e^\epsilon
d(\Omega)$. But this means that $d(\overline{\Omega})\leq
e^\epsilon d(\Omega)$, and since $\epsilon >0$ was arbitrary we
conclude that $d(\overline{\Omega})=d(\Omega)$. {\bf q.e.d.}

\section{References}

\begin{enumerate}
\item[{[B]}] A. Banyaga,  {\sl Sur la structure du groupe des
diffeomorphismes qui pr\'e\-servent une forme symplectique}, Comm. Math. Helv. {\bf 53} (1978), 174-227.
\item[{[C]}] Y. Chekanov, {\sl Lagrangian intersections, symplectic
energy, and area of holomorphic curves}, Duke Math. J. 95 (1998), 213-226.
\item[{[EH]}] I. Ekeland, H. Hofer, {\sl Symplectic topology and
Hamiltonian dynamics}, Math. Zeit. 20 (1990), 355-378.
\item[{[E]}] Y. Eliashberg, {\sl Filling by holomorphic discs
and its applications}, London Math. Society Lecture Notes,
Series 151 (1991), 45-67.
\item[{[G]}] M. Gromov, {\sl Pseudo-holomorphic curves in
symplectic manifolds}, Invent. Math. 82 (1985), 307-347.
\item[{[He]}] D. Hermann, {\sl Non-equivalence of symplectic capacities
for open sets with restricted contact type boundary},
Preprint, Orsay 1998.
\item[{[Ho]}] H. Hofer, {\sl On the topological properties
of symplectic maps}, Proc. R. Soc. Edinb., Sect. A, 115 (1990), 25-38.
\item[{[HZ]}] H. Hofer, E. Zehnder, {\it Symplectic invariants
and Hamiltonian dynamics}, Birkh\"auser, Basel 1994.
\item[{[LM]}] F. Lalonde, D. McDuff, {\sl The geometry of symplectic energy},
Ann. Math. {\bf 141} (1995), 349-371.
\item[{[MD1]}] D.~McDuff, {\sl The local behaviour of holomorphic
curves in almost complex 4-manifolds,} J. Diff. Geom.
34 (1991), 143-164.
\item[{[MD2]}] D.~McDuff, {\sl Symplectic manifolds
with contact type boundaries}, Invent. Math 103 (1991), 651-671. 
\item[{[MS]}] D.~McDuff, D. Salamon, {\it Introduction to symplectic topology},
Oxford Univ. Press, Oxford (1995).
\item[{[MDT]}] D.~McDuff, L. Traynor, {\sl The 4-dimensional
symplectic camel and related results}, in Symplectic
geometry (ed. D. Salamon), London Math. Society Lecture Notes 192,
169-182, Cambridge University Press 1993.
\item[{[O]}] Y.-G. Oh, {\sl Removal of boundary singularities
of pseudo-holomophic curves with Lagrangian boundary conditions,}
Comm. Pure Appl. Math. 45 (1992), 121-139.
\item[{[P]}] L. Polterovich, {\it The geometry of the group of
symplectic diffeomorphisms}, Birkh\"auser, Basel 2001.
\item[{[V]}] C. Viterbo, {\sl Metric and isoperimetric
problems in symplectic geometry}, J. AMS. 13 (2000), 411-431.
\end{enumerate}

\end{document}